\documentclass[11pt]{amsart}
\usepackage{amsmath,amssymb,amsthm}
\usepackage{graphicx}
\usepackage{hyperref}
\hypersetup{colorlinks=true,citecolor=blue,linkcolor=blue}
\newtheorem{theorem}{Theorem}[section]
\newtheorem{lemma}[theorem]{Lemma}

\begin{document}
\title{Existence of Classical
Solution of the  Boussinesq  Equation}
\author{Wu  Shu-hong }
\address{Department of Mathematics,School of Mathematics and  Statistics,\\ Wuhan
University of Technology,  Wuhan  P.R.C.,430070.   China .}
\email{wutwsh@163.com}
\date{}\begin{abstract}We   make use of compression  fixed point principle to discuss    the  existence of classical
 Solution of  the Boussinesq
 equation  to dimensions  $n\geq  8$   and   the  existence     of classical
 Solution of  the  Navier-Stokes
 equations  on   ${ \mathbf{R}}
^n$  to the dimensions  $n\geq  4$ and
the  existence    of classical
 Solution of  Navier-Stokes
 equations  on   $ {\mathbf{R}}
^n/{\mathbf{Z}}^n
$  to the dimensions  $n\geq  3$.Using  this  method,  if  we  know 
 a  suitable initial value ,  we  can  calculate the limit solution.
\end{abstract}
\keywords{compression  fixed point principle;
   Boussinesq
 equations;     Navier-Stokes equations; classical solution, existence}
 \tableofcontents
\maketitle
\newpage
\section{Introduction}
The $n$-dimensional incompressible Boussinesq equations  and Navier-Stokes equations
     concerned here can be represented in the form
      \begin{equation}\left \{ \begin{array}{lllll} u_  t + u \cdot
       \nabla u- \nu\triangle u+ \nabla p=A\rho  +f, & t>0, x\in {\bf{R}}^n \\
\rho _t+ u\cdot \nabla  \rho -   \kappa  \triangle  \rho=g,\\
 \textrm{div} u=0,\\
    \left . u(x,t)\right |_{t=0}=\varphi(x),     \left . \rho(x,t)
    \right |_{t=0}=\eta(x), p(x_0,t)=p_0(t),
      \end{array}             \right.  \end{equation}
        \begin{equation}\left \{ \begin{array}{lllll} u_  t + u \cdot
       \nabla u- \nu\triangle u+ \nabla p= f, & t>0, x\in {\bf{R}}^n \\
 \textrm{div} u=0,\\
    \left . u(x,t)\right |_{t=0}=\varphi(x),    p(x_0,t)=p_0(t),
      \end{array}             \right.  \end{equation}
respectively,where    $\nu,\kappa >0,u,\rho: {\bf{R}}^n\times {\bf{R}}_+
\rightarrow
 {\bf{R}}^n$  are vector fields; $p$ is
 scalar,  $A$  is   $n\times n$    matrix.It  is  easy   to see 
  that   Navier-Stokes equations  is a special case of Boussinesq equations  as   $\rho(x,t)=0,g(x)=\eta(x)=0.$

      The Boussinesq system   describes  the  influence of  convection
       (or convection-diffusion)
             phenomenon in a  viscous or inviscid  fluid. The Boussinesq
Equations are extensively used
in the atmospheric sciences and oceanographic turbulence in which rotation and
 stratification
are important (see, for example, [1]) and references therein). Thus, over the past
 few years,
the Boussinesq Equations were extensively studied  theoretically. Existence of
 global  solutions  is a challenging open problem ([2]-[9]).

  In this paper,  we     make use of compression  fixed point principle
 to discuss    the  existence of classic
 Solution of  the Boussinesq
 equation  to dimensions  $n\geq  8$   and   the  existence     of classic
 Solution of  the  Navier-Stokes
 equations  on   ${ \mathbf{R}}
^n$  to the dimensions  $n\geq  4$ and
the  existence    of classic
 Solution of  Navier-Stokes
 equations  on   $ {\mathbf{R}}
^n
/{
\mathbf{Z}}
^n
$  to the dimensions  $n\geq  3$.Using  this  method, 
 if  we  know  a  suitable initial value ,  we  can  calculate the limit solution.

\section{Some Lemmas}

 \begin{lemma} Let   $g,h\in  K,\eta,\varphi\in J, \psi(x,t)=\sum\limits_{m\in {\bf{Z}}^n} \widehat{\psi}(m,t)  \exp
 (\sum\limits_{j=1}^n  \frac{2im_j\pi x_j}{l_j}),$     $
\phi(x)=\sum\limits_{m\in {\bf{Z}}^n} \widehat{\phi}(m)  \exp
 (\sum\limits_{j=1}^n \frac{2im_j\pi x_j}{l_j}),\overline{I}(\phi)= \sum\limits_{ m\in{\bf{Z}}^n }
  |\widehat{\phi}(m)| ,  I(\psi)
= \sum\limits_{ m\in{\bf{Z}}^n } \int_0^\infty |\widehat{\psi}(m,t)|dt,$
 $$K=\{\alpha\in  C^{\infty}({\mathbf{R}^n\times{\mathbf{R}^+}})|\alpha (x,t)=
 \sum\limits_{\beta\in  {\mathbf{Z}^n}} \widehat{\alpha}(\beta,t)\exp(\sum\limits_{j=1}^n\frac{2\beta_j\pi 
  x_ji}{l_j}),I(\alpha)<+\infty\},$$
 $$J=\{\alpha\in  C^{\infty}({\mathbf{R}^n})|\alpha (x)=\sum\limits_{\beta\in  {\mathbf{Z}^n}}
  \widehat{\alpha}(\beta)\exp(\sum\limits_{j=1}^n\frac{2\beta_j\pi  x_ji}{l_j}),\overline{I}(\alpha)<+\infty\},$$
 \begin{eqnarray*} &\/&\rho(x,t)= \sum\limits_{m\in{\bf{Z}}^n} \exp
 [\sum\limits_{j=1}^n [ \frac{2im_j\pi x_j}{l_j} -
 (\frac{2m_j\pi}{l_j})^2\kappa t]  \left[\widehat{\eta} (m)    +
  \int_0^t  \widehat{g}(m,\tau)  \exp [\sum\limits_{j=1}^n
 (\frac{2m_j\pi}{l_j})^2\kappa \tau]d \tau\right]  \\&\/&
  u(x,t)= \sum\limits_{m\in{\bf{Z}}^n} \exp
 [\sum\limits_{j=1}^n [ \frac{2im_j\pi x_j}{l_j} -
 (\frac{2m_j\pi}{l_j})^2\nu t]  \left[\widehat{\varphi}(m)    +
  \int_0^t  \widehat{h}(m,\tau)  \exp [\sum\limits_{j=1}^n
\tilde{}  (\frac{2m_j\pi}{l_j})^2\nu \tau]d \tau\right] . \end{eqnarray*}  Then
\begin{eqnarray*} &\/&I(u\frac{\partial \rho}{\partial x_k})
  \leq \frac{l_k}{2\pi  \kappa}  \sum\limits_{p\in{\bf{Z}}^n}  
    \left[|\widehat{\varphi}(p)|    +
  \int_0^\infty |\widehat{h}(p,\tau)|  d \tau\right]\cdot  
  \sum\limits_{q\in{\bf{Z}}^n}\left[|\widehat{\eta}(q) |   +
  \int_0^\infty  |\widehat{g}(q,\tau) |  d \tau\right]\\  
   &=&\frac{l_k}{2\pi  \kappa}  [ \overline{I}(\varphi)  
    +I(h)] [ \overline{I}(\eta)   +I(g)].\end{eqnarray*}
   \end{lemma}

   Proof:\hspace{0.3cm}  As
   \begin{eqnarray*} &\/&\int_0^\infty  dz\int_0^z dx\int_0^z F(x,y,z)dy=
    \int_0^\infty  dx\int_x^\infty dz\int_0^z F(x,y,z)dy\\&=&  \int_0^\infty
  dx\int_0^\infty dy\int_{\max\{x,y\}}^\infty F(x,y,z)dz , \end{eqnarray*} 
  \begin{eqnarray*} &\/&\widehat{u\frac{\partial \rho}{\partial x_k}}(m,t)\\&=& 
    \sum\limits_{p,q\in{\bf{Z}}^n,m=p+q,q_k\neq 0} \left\{ \widehat{\varphi}(p  )
     \exp [-\sum\limits_{j=1}^n
(\frac{2p_j\pi}{l_j})^2\nu t]  +
  \int_0^t  \widehat{h}(p,\tau)  \exp [\sum\limits_{j=1}^n
(\frac{2p_j\pi}{l_j})^2\nu(\tau-t)  ]d \tau\right\} \\  &\/& 
   \cdot \frac{2iq_k\pi }{l_k} \left\{\widehat{\eta}(q) \exp
 \{-\sum\limits_{j=1}^n
 (\frac{2q_j\pi}{l_j})^2\kappa t    +
  \int_0^t  \widehat{g}(q,\tau)  \exp [\sum\limits_{j=1}^n
 (\frac{2q_j\pi}{l_j})^2\kappa  ( \tau-t)]\}d \tau \right\}\\
&\/&\int_0^\infty  \left|\widehat{u\frac{\partial \rho}{\partial x_k}}(m,t)
\right|dt\leq \sum\limits_{p,q\in{\bf{Z}}^n,m=p+q,q_k\neq 0}\frac{2|q_k|\pi }{l_k}\int_0^\infty\exp
 [\sum\limits_{j=1}^n [  -
 (\frac{2p_j\pi}{l_j})^2\nu t -
 (\frac{2q_j\pi}{l_j})^2\kappa t]\\  &\/&    \left[|\widehat{\varphi}(p)|    +
  \int_0^t  |\widehat{h}(p,\tau)|  \exp [\sum\limits_{j=1}^n
(\frac{2p_j\pi}{l_j})^2\nu \tau]d \tau\right]\cdot \left[|\widehat{\eta}(q) |   +
  \int_0^t  |\widehat{g}(q,\tau) | \exp [\sum\limits_{j=1}^n
 (\frac{2q_j\pi}{l_j})^2\kappa \tau]d \tau\right] dt\\
&\leq&  \sum\limits_{p,q\in{\bf{Z}}^n,m=p+q,q_k\neq 0} \frac{2|q_k|\pi }{l_k} 
 \left[ \sum\limits_{j=1}^n[
 (\frac{2p_j\pi}{l_j})^2\nu +
 (\frac{2q_j\pi}{l_j})^2\kappa] \right]^{-1}\\  &\/&    \left[|\widehat{\varphi}(p)|    +
  \int_0^\infty |\widehat{h}(p,\tau)|  d \tau\right]\cdot \left[|\widehat{\eta}(q) |   +
  \int_0^\infty  |\widehat{g}(q,\tau) |  d \tau\right] \\
&\leq&   \sum\limits_{p,q\in{\bf{Z}}^n,m=p+q,q_k\neq 0} \frac{l_k}{2\pi  \kappa}   
  \left[|\widehat{\varphi}(p)|    +
  \int_0^\infty |\widehat{h}(p,\tau)|  d \tau\right]\cdot \left[|\widehat{\eta}(q) |   +
  \int_0^\infty  |\widehat{g}(q,\tau) |  d \tau\right]\end{eqnarray*}
\begin{eqnarray*}&\/&I(u\frac{\partial \rho}{\partial x_k})\leq \frac{l_k}{2\pi  \kappa} 
   \sum\limits_{p\in{\bf{Z}}^n}    \left[|\widehat{\varphi}(p)|    +
  \int_0^\infty |\widehat{h}(p,\tau)|  d \tau\right]\cdot  
  \sum\limits_{q\in{\bf{Z}}^n}\left[|\widehat{\eta}(q) |   +
  \int_0^\infty  |\widehat{g}(q,\tau) |  d \tau\right]\\ 
    &\leq &\frac{l_k}{2\pi  \kappa}  [ \overline{I}(\varphi)   +I(h)] 
    [ \overline{I}(\eta)+I(g)].\end{eqnarray*}\\

  \begin{lemma}  (1) If $\mu>0$,  the solution
  of the system
\begin{eqnarray*}\left \{\begin{array}{lll}{ \frac{\partial u }
{\partial t}}-\mu\triangle
u=f(x,t),x\in {\bf{R}}^n,t>0,\\
 u(x,0)=\varphi(x)
 \end{array}  \right.\end{eqnarray*} is $$u(x,t)=
 \int_{ {\bf{R}}^n}\frac{\varphi(y )e^{\frac{-|x-y|^2}{4 \mu t}}}
 {(4\pi \mu t)^{ \frac{n}{2}}} dy+ \int_0^t\int_{ {\bf{R}}^n}
 \frac{f(y,\tau)e^{\frac{-|x-y|^2}
 {4 \mu(t-\tau)}}}{[4\pi \mu (t-\tau)]^{ \frac{n}{2}}}d\tau dy,$$
(2)If   $n\geq 3,  g\in   C^2({\mathbf{R}^n})$, there  exist $C>0,\alpha>n$   such  that   
  $|g(y)|\leq   C(1+|y|)^{-\alpha},  \forall  y\in  {\mathbf{R}^n},$   $$u(x)=u_0-\int_{ {\bf{R}}^n}\frac{\Gamma
 ( \frac{n}{2})(|x-y|^{2-n}-|x_0-y|^{2-n})g(y )} { 2(2-n)\pi^{ \frac{n}{2}}} dy$$
 is  the solution of the system
$$\triangle
u (x)=g(x),x\in {\bf{R}}^n, u(x_0)=u_0 $$  and   $u\in  C^\infty  ({\bf{R}}^n).$
\begin{eqnarray*}(3) &\/&u(x,t)=
 \sum\limits_{m\in{\bf{Z}}^n}  \exp [\sum\limits_{j=1}^n (
  \frac{2im_j\pi x_j}{l_j} -
 (\frac{2m_j\pi}{l_j})^2\mu t)] \left[  \widehat{ \varphi}({\mathbf{m}}) +
  \int_0^t\widehat{f}({\mathbf{m}},\tau) \exp (\sum\limits_{j=1}^n
 (\frac{2m_j\pi}{l_j})^2\mu \tau)  d \tau  \right]  \end{eqnarray*}
   is periodic solution (period: $l_k$) of $x_k(1\leq k\leq n)$ of system
    $$\left \{\begin{array}{lll}{ \frac{\partial u }{\partial t}}-\mu\triangle
u=f(x,t),x\in  \prod\limits_{i=1}^n [0,l_i],t>0,\\
 u(x,0)=\varphi(x).
 \end{array}  \right.$$  \begin{eqnarray*}(4)\hspace{4cm}
   p(x)= p_0-\sum\limits_{m\in{\bf{Z}}^n/ \{\mathbf{0}\} }
  \frac{ \widehat{g}_m [\exp  (\sum\limits_{j=1}^n  \frac{2im_j\pi x_j}{l_j})-
  \exp  (\sum\limits_{j=1}^n  \frac{2im_j\pi x_{0,j}}{l_j})] }{4\pi^2\sum\limits_{j=1}^n
     \frac{m_j^2}{l_j^2} }     .  \end{eqnarray*}
   is periodic solution (period: $l_k$) of $x_k(1\leq k\leq n)$ of system
   $$ \triangle
p (x)=g(x),x\in  \prod\limits_{i=1}^n [0,l_i], p(x_0)=p_0.$$\end{lemma}

    Proof:\hspace{0.3cm}    Proof:\hspace{0.3cm}  (1),(3),(4)  are  very simple,
     (2) can be proved just as Theorem 4.1 in [10].\\

\section{Main Result} 

 \begin{theorem}  Let   $   1\leq k\leq n,j=0,1,
 \nu,\kappa,\lambda,D>0, \varphi_k, \eta_k,f_k,g_k,r_k,h_k\in H(\lambda,D),\zeta=(r_1,\cdots,r_n,h_1,
 \cdots,h_n)^\top,
    \varphi=(\varphi_1,\cdots, \varphi_n)^\top, \textrm{div} \varphi =0,|f_k|_{j+}=\int_{ {\bf{R}}^n}
      |\xi|^j \sup\limits_{t\in   {\bf{R}}^+}
          | \widehat{f_k}(\xi,t ) |    d\xi,|\widehat{f_k}|_{j-}=\sup\limits_{\xi\in  {\bf{R}}^n,t\in 
           { {\bf{R}}}^+}  |\xi|^j\cdot|\widehat{f}_k(\xi,t)|,  |h|_{j+}$   $= \max\limits_{ 1\leq k\leq n}
            | h_k|_{j+},  |h|_{j-}= \max\limits_{ 1\leq k\leq n} | h_k|_{j-},||h||= \max\{ \lambda^{-1} 
            | h|_{0+},|h|_{0-}\},$   $  ||\zeta|| =   || h||+|| r ||  ,  |\varphi_k|_{1-}<+\infty, 
             |\varphi_k|_{1+}<+\infty, |\eta_k|_{1-}<+\infty, |\eta_k|_{1+}<+\infty, Q_k(\omega,t)
             =\widehat{r}_k(\omega,t)+\widehat{f}_k(\omega,t)+e_kA\{\widehat{\eta}(\omega )e^{ -\kappa t|
   \omega|^2}+\int_0^t [\widehat{g}(\omega,\tau)+\widehat{h}(\omega,\tau)]$    $
   e^{ -\kappa(t-\tau)|\omega|^2 }d\tau \},$   $P_k(\omega,t)=Q_k(\omega,t) - 
   \sum\limits_{m=1}^n\frac{\omega_m\omega_k}{|\omega|^2}$   $ Q_m(\omega,t),S_n= 2\pi  
   \int_0^\pi\sin ^{n-2}\varphi_1  d\varphi_1 \int_0^\pi \sin ^{n-3}\varphi_2 d\varphi_2\cdots  
    \int_0^\pi \sin \varphi_{n-2} d\varphi_{n-2}$   $ =\frac{\pi^{\frac{n}{2}}n}{\Gamma(\frac{n}{2}+1)},$
 $$     H(\lambda,D)= \{\alpha\in  H|  |\alpha|_{0+}\leq \lambda    D, |\alpha|_{0-} \leq D,  
  |\alpha  (x,t)|\leq C(1+|x|)^{-c},C>0,c>n\},$$ \begin{eqnarray*}&\/&M_1(\lambda,D)\\&=&
| \eta|_{1-}[ | \varphi|_{0+} +\nu^{-1}(n+1)(D+||f||+nB||\eta||)
(\frac{ S_n  }{ n-2 }+\lambda)\\&\/&+(n+1)nB\nu^{-1}\kappa^{-1}(D+||g|| )( \frac{S_n  }{n-4} +\lambda)]+  
\kappa^{-1} |\varphi|_{0-}(D+||g|| )(  \frac{S_n  }{n-1} +\lambda) \\  &\/& +  \nu^{-1}\kappa^{-1} (n+1) (D+||f||  +nB||
 \eta||)(D+||g|| )(\frac{S_n  }{n-3} +\lambda)  \\&\/&
 +n (n+1) B\kappa^{-2}\nu^{-1}
   ( ||g|| +D)^2\cdot  (\frac{   S_n}{n-5}+\lambda) \\ &\/&M_2(\lambda,D)\\&=& | \varphi|_{1-}[|\varphi|_{0+} +(n+1)\nu^{-1}
(D+||f|| +nB|| \eta||)( \frac{S_n  }{ n-2 } +\lambda)\\&\/&+n(n+1)B\nu^{-1}\kappa^{-1}(D+||g|| )
( \frac{S_n  }{n-4} +\lambda)]+(n+1)\nu^{-1}| \varphi|_{0-}(D+||f|| +nB|| \eta||)(\frac{S_n  }{n-1}+\lambda) 
\\&\/&+(n+1)^2\nu^{-2}(||f||   +D+nB| |\eta||_{0})^2 (\frac{  S_n}{n-3}+ \lambda) \\  &\/& + 2 n(n+1)^2B \nu^{-2}\kappa^{-1} 
 (D+||f||_{0}    +nB ||\eta||_{0})   (D+||g||_{0} )
 (\frac{  S_n }{n-5}+ \lambda)\\&\/&+ n(n+1)|\varphi|_{0-}B\nu^{-1}\kappa^{-1}(D+||g|| )( \frac{S_n  }{n-1} +\lambda) ]
 \\&\/& + n^2 (n+1)^2B^2\nu^{-2}\kappa^{-2}
  (||g|| +D)^2  (\frac{  S_n}{n-7}+\lambda )  ,\end{eqnarray*}  \begin{eqnarray*}&\/&N_1(\lambda,D)\\&=& [ |\eta|_{1+} +
     \kappa^{-1}(||g||  +D)( \frac{ S_n}{ n-1 }+\lambda)]
[ | \varphi|_{0+}  + (n+1)\nu^{-1}(D+ ||f|| + nB||\eta|| ) (  \frac{  S_n}{n-2}+\lambda)\\&\/&+
n(n+1)B\kappa^{-1}\nu^{-1}(||g||+D)(\frac{S_n}{n-4}+\lambda)]\\  &\/& N_2(\lambda,D)\\ &=&[
   | \varphi|_{0+}  + (n+1)\nu^{-1}(D+ ||f|| + nB||\eta|| ) (  \frac{  S_n}{n-2}+\lambda)+
n(n+1)B\kappa^{-1}\nu^{-1}(||g||+D)(\frac{S_n}{n-4}+\lambda)]\\&\/& [ | \varphi|_{1+}  +
 (n+1)\nu^{-1}(D+ ||f|| + nB||\eta|| ) (  \frac{  S_n}{n-1}+\lambda) +
n(n+1)B\kappa^{-1}\nu^{-1}(||g||+D)(\frac{S_n}{n-3}+\lambda)] ,\end{eqnarray*}\begin{eqnarray*} 
  M(\lambda,D)=(2\pi)^{-n}n \max \{  M_1(\lambda,D),  M_2(\lambda,D)\},
 N(\lambda,D)&=&n (2\pi)^{-n}\max \{  N_1(\lambda,D), N_2(\lambda,D)\}, \end{eqnarray*}
 \begin{eqnarray*}
    P_{1-}^r  =  (n+1)\nu^{-1} [
| \eta|_{1-}      (\frac{ S_n  }{ n-2 }+\lambda)+ \kappa^{-1}(D+||g||)   (\frac{S_n  }{n-3} +\lambda)] 
 \\  P_{1-}^h  =  (n+1)nB\nu^{-1}\kappa^{-1}[| \eta|_{1-}( \frac{S_n  }{n-4} +\lambda)+  \kappa^{-1}( ||g||+D)
 (  \frac{S_n  }{n-5} +\lambda) ]  , \end{eqnarray*}
    \begin{eqnarray*}
  &\/&  Q_{1-}^h =
  [\kappa^{-1} |\varphi|_{0-}(  \frac{S_n  }{n-1} +\lambda) +  \kappa^{-1}  \nu^{-1}(n+1) (D+||f|| 
   +nB|| \eta||)(\frac{S_n  }{n-3} +\lambda)  \\&\/&
 +n (n+1) B\kappa^{-2}\nu^{-1}
   ( ||g|| +D)\cdot  (\frac{   S_n}{n-5}+\lambda)] , \end{eqnarray*}
\begin{eqnarray*}&\/&
    P_{2-}^r =  (n+1)[| \varphi|_{1-} \nu^{-1}(\frac{S_n  }{n-2}+\lambda)+ (n+1)\nu^{-2}\kappa^{-1}(D+||g||)
    ( \frac{S_n  }{n-5} +\lambda) \\&\/&+ (n+1) \nu^{-2} (D+||f||    +nB ||\eta||)
 (\frac{  S_n }{n-3}+ \lambda)] \\&\/&P_{2-}^h = (n+1)\kappa  ^{-1} \nu^{-1}[| \varphi|_{1-} (\frac{S_n  }
 {n-4}+\lambda)+ (n+1)\nu^{-1}\kappa^{-1}(D+||g||)( \frac{S_n  }{n-7} +\lambda)
  \\&\/&+ (n+1) \nu^{-1} (D+||f||    +nB ||\eta||)
 (\frac{  S_n }{n-5}+ \lambda)]
 ,\end{eqnarray*}
    \begin{eqnarray*}&\/&
    Q_{2-}^r =(n+1)\nu^{-1}[| \varphi|_{1-} (\frac{S_n  }{n-2}+\lambda)+ (n+1)\nu^{-1}\kappa^{-1}(D+||g||)
    ( \frac{S_n  }{n-5} +\lambda) \\&\/&+ (n+1) \nu^{-1} (D+||f||    +nB ||\eta||)
 (\frac{  S_n }{n-3}+ \lambda)]  \\&\/&
    Q_{2-}^h =  (n+1)\kappa  ^{-1}\nu^{-1}[| \varphi|_{1-} (\frac{S_n  }{n-4}+\lambda)+ (n+1)\nu^{-1}
    \kappa^{-1}(D+||g||)( \frac{S_n  }{n-7} +\lambda) \\&\/&+ (n+1) \nu^{-1} (D+||f||    +nB ||\eta||)
 (\frac{  S_n }{n-5}+ \lambda)]\\,\end{eqnarray*}
   \begin{eqnarray*}&\/&P_{1+}^r =(n+1)\nu^{-2}[ |\eta|_{1+}( \frac{ S_n}{ n-2 }+\lambda) + 
    \kappa^{-1}(||g||  +D)( \frac{ S_n}{ n-1 }+\lambda)]( \frac{ S_n}{ n-2 }+\lambda)\\&\/&P_{1+}^h 
    =(n+1)\nu^{-2}\kappa^{-1}[ |\eta|_{1+}( \frac{ S_n}{ n-2 }+\lambda) + \kappa^{-1}(||g||  +D)
    ( \frac{ S_n}{ n-1 }+\lambda)]( \frac{ S_n}{ n-4 }+\lambda)  ,  \end{eqnarray*}
  \begin{eqnarray*}&\/&Q_{1+}^h=\kappa^{-1}
[ | \varphi|_{0+}  + (n+1)\nu^{-1}(D+ ||f|| + nB||\eta|| ) (  \frac{  S_n}{n-2}+\lambda)\\   
 &\/&+\kappa^{-1}\nu^{-1}
n(n+1)B(||g||+D)(\frac{S_n}{n-4}+\lambda)] (\frac{S_n}{n-1}+\lambda)   ,  \end{eqnarray*}
 \begin{eqnarray*}&\/&P_{2+}^r =    [ | \varphi|_{1+}  + (n+1)\nu^{-1}(D+ ||f|| + nB||\eta|| ) 
  (  \frac{  S_n}{n-1}+\lambda) \\&\/&+\kappa^{-1}\nu^{-1}
n(n+1)B(||g||+D)(\frac{S_n}{n-3}+\lambda)] (n+1)\nu^{-1}
(\frac{S_n}{n-2}+\lambda). \end{eqnarray*} \begin{eqnarray*}&\/&P_{2+}^h =   
   [ | \varphi|_{1+}  + (n+1)\nu^{-1}(D+ ||f|| + nB||\eta|| ) (  \frac{  S_n}{n-1}+\lambda) 
   \\&\/&+\kappa^{-1}\nu^{-1}
n(n+1)B(||g||+D)(\frac{S_n}{n-3}+\lambda)] (n+1)\nu^{-1}
\kappa^{-1}nB(\frac{S_n}{n-4}+\lambda). \end{eqnarray*}\begin{eqnarray*}&\/&Q_{2+}^r= 
  [ | \varphi|_{0+}  + (n+1)\nu^{-1}(D+ ||f|| + nB||\eta|| ) (  \frac{  S_n}{n-2}+\lambda)
  \\&\/&+\kappa^{-1}\nu^{-1}
n(n+1)B(||g||+D)(\frac{S_n}{n-4}+\lambda)] (n+1)\nu^{-1} (  \frac{  S_n}{n-3}+\lambda) . 
\end{eqnarray*}\begin{eqnarray*}&\/&Q_{2+}^r= [ | \varphi|_{0+}  + (n+1)\nu^{-1}(D+ ||f|| + nB||\eta|| )
   (  \frac{  S_n}{n-2}+\lambda)\\&\/&+\kappa^{-1}\nu^{-1}
n(n+1)B(||g||+D)(\frac{S_n}{n-4}+\lambda)] (n+1)\nu^{-1} \kappa^{-1}nB (  \frac{  S_n}{n-1}+\lambda) . 
\end{eqnarray*}$$\varepsilon=n( P_{2-}^r+P_{2-}^h+Q_{2-}^r+Q_{2-}^h+P_{1-}^r+P_{1-}^h+Q_{1-}^h+
P_{2+}^r+P_{2+}^h+Q_{2+}^r+Q_{2+}^h+
P_{1+}^r+P_{1+}^h+Q_{1+}^h).$$If  there  exist  $\lambda,D>0$  such  that $M(\lambda,D)\leq  D,N(\lambda,D)\leq 
 \lambda^{-1} D,\varepsilon  <1,$
    then there  exist  solutions   $ u_k,\rho_k \in  H(1\leq k\leq n),p\in  C^\infty({\bf{R}}^n),\sup\limits_{  t\geq 0}
 \int_{{
 \mathbf{R}}^n} |u_k(x,t)|^2
  dx<+\infty $   for system (1)   as $n\geq 8$ .

   If there  exist $\lambda,D>0$   such  that $   
   \varepsilon=n\max\{\varepsilon_1,\varepsilon_2\}<1,$ \begin{eqnarray*}\varepsilon_1&=&  (2\pi)^{n}
(n+1) \nu^{-1}[(  \frac{  S_n}{ n-2 }+  \lambda)| \varphi|_{1-} + ( \frac{  S_n}{ n-1 }+  
\lambda)| \varphi|_{0-} +2\nu^{-1}(n+1)(D+||f||  )(  \frac{  S_n}{ n-3 }+  \lambda)], \end{eqnarray*}
 \begin{eqnarray*}\varepsilon_2&=& (2\pi)^{n}  \nu^{-1}(n+1)\lambda^{-1}(  \frac{  S_n}{ n-2 }+  
  \lambda)[ | \varphi|_{1+}+\nu^{-1}(n+1)(  \frac{  S_n}{ n-1 }+  \lambda) (D+||f||) ]+\\ &\/&
 \nu^{-1}(n+1)\lambda  ^{-1}(  \frac{  S_n}{ n-1 }+  \lambda)[| \varphi|_{0+} +\nu^{-1}(n+1)(D+||f||  )
 (  \frac{  S_n}{ n-2 }+  \lambda)],\end{eqnarray*}
 \begin{eqnarray*}&\/&[ | \varphi|_{0+}  + \nu^{-1}(n+1)(D+||f||  )(  \frac{  S_n}{ n-2 }+ 
   \lambda) ][ |\varphi|_{1+}  + \nu^{-1}(n+1)(D+||f||  )(  \frac{  S_n}{ n-1 }+  \lambda)]\leq 
   \frac{(2\pi)^n\lambda D}{n}, \end{eqnarray*}\begin{eqnarray*}&\/&
| \varphi|_{1-}|\varphi|_{0+}+\nu^{-1}(n+1)| \varphi|_{1-}(D+||f||  )(  \frac{  S_n}{ n-2 }+ 
 \lambda) \\ &\/& +\nu^{-1}(n+1)| \varphi|_{0-}(D+||f||  )(  \frac{  S_n}{ n-1 }+  \lambda) +\nu^{-2}(n+1)^2( 
   \frac{  S_n}{ n-3 }+  \lambda)(||f||   +D)^2 \leq \frac{(2\pi)^n D}{n}  ,\end{eqnarray*}
Then there  exist  solutions   $ u_k\in H (1\leq k\leq n),p\in   C^\infty({\bf{R}}^n\times {\mathbf{R}}^+),
\sup\limits_{  t\geq 0}
 \int_{{{
 \mathbf{R}}}^n} |u_k(x,t)|^2
  dx<+\infty $   for system (1)   as $n\geq 4$ . \end{theorem}

 Proof:\hspace{0.3cm} Let   $ \lambda,D>0,H(\lambda,D)= \{\alpha\in  H|  |\alpha|_{0+}\leq \lambda  ^{-1}D, 
 |\alpha|_{0-} \leq D\}, \zeta=(h_1,\cdots, h_n, r_1,\cdots,
               r_n)^\top   $   $\in  H^{2n}(\lambda,D),h=(h_1,\cdots, h_n )^\top,r=(  r_1,\cdots,
               r_n)^\top,g=(g_1,\cdots, g_n )^\top,f=(f_1,\cdots, f_n )^\top, $ $             
                 A=(a_{i,j}),B=\max\limits_{1\leq i,j\leq n}|a_{i,j}|,$
    $$ \rho_k(x,t)=\int_{ {\bf{R}}^n}\frac{\eta_k(y
  )e^{\frac{-|x-y|^2}{4\kappa   t}}}{(4\pi  \kappa t)^{ \frac{n}{2}}}
   dy+\int_0^t\int_{ {\bf{R}}^n} \frac{[g_k(y,\tau)+h_k(y,\tau)]e^{\frac{-
   |x-y|^2}{4 \kappa (t-\tau)}}}{[4\pi  \kappa (t-\tau)]^{ \frac{n}{2}}}d\tau dy    ,
  $$  \begin{eqnarray*}&\/&p(x,t)=p_0(t)-\int_{ {\bf{R}}^n}\frac{\Gamma
 ( \frac{n}{2})(|x-y|^{2-n}-|x_0-y|^{2-n})[\textrm{div}r(y,t ) + \textrm{div}f(y,t )+
\sum\limits_{i=1}^n\frac{\partial [e_i  A \rho(y,t )]}{\partial  y_i}]} { 2(2-n)\pi^{ \frac{n}{2}}} dy 
 \\  &=&p_0(t)-\int_{ {\bf{R}}^n}\frac{\Gamma
 ( \frac{n}{2})(|x-y|^{2-n}-|x_0-y|^{2-n})[\textrm{div}r(y,t ) + \textrm{div}f(y,t )+
\sum\limits_{i=1}^n\sum\limits_{j=1}^na_{i,j}\frac{\partial  \rho_j(y,t )}{\partial  y_i}]} { 2(2-n)\pi^{ \frac{n}{2}}} dy 
   , \end{eqnarray*} $$ u_k(x,t)=
  \int_{ {\bf{R}}^n}\frac{\varphi_k(y )e^{\frac{-|x-y|^2}{4 \nu t}}}{(4\pi \nu t)
  ^{ \frac{n}{2}}} dy+\int_0^t\int_{ {\bf{R}}^n} \frac{[r_k(y,\tau)+f_k(y,\tau)+
  e_k A \rho(y,\tau ) -\frac{\partial p(y,\tau)}{\partial y_k}]e^{\frac{-|x-y|^2}
  {4 \nu(t-\tau)}}}{[4\pi \nu (t-\tau)]^{ \frac{n}{2}}}\acute{}d\tau dy   .$$
$ Q_k(\omega,t)=\widehat{r}_k(\omega,t)+\widehat{f}_k(\omega,t)+e_kA[\widehat{\eta}(\omega )e^{ -\kappa t|
   \omega|^2}+\int_0^t [\widehat{g}(\omega,\tau)+\widehat{h}(\omega,\tau)]
   e^{ -\kappa(t-\tau)|\omega|^2 }d\tau ],P_k(\omega,t)=Q_k(\omega,t) - 
   \sum\limits_{m=1}^n\frac{\omega_m\omega_k}{|\omega|^2}Q_m(\omega,t).$
 Then
 $$\triangle p(x,t)=\textrm{div}r(x,t ) + \textrm{div}f(x,t )+\sum\limits_{i=1}^n\sum\limits_{j=1}^na_{i,j}
 \frac{\partial  \rho_j(x,t )}{\partial  x_i},p(x_0,t)=p_0(t), $$  $$ \widehat{\rho}_k(\omega,t)= 
 \widehat{\eta}_k(\omega )e^{ -\kappa t|
   \omega|^2}+\int_0^t [\widehat{g}_k(\omega,\tau)+\widehat{h}_k(\omega,\tau)]
   e^{ -\kappa(t-\tau)|\omega|^2 }d\tau ,$$$$ \widehat{p}(\omega,t)=-\sum\limits_{m=1}^n\frac{i\omega_m}
   {  |\omega |^2 }
       [ \widehat{r_m}(\omega,t)+ \widehat{f_m}(\omega,t)+e _m A  \widehat{\rho}(\omega,t)], $$ $$     
         \widehat{\frac{\partial p }{\partial x_k}}
 (\omega,t)=\sum\limits_{m=1}^n\frac{\omega_m\omega_k}{|\omega|^2}
        [  \widehat{r_m}(\omega,t)+ \widehat{f_m}(\omega,t)+e _m A  \widehat{\rho}(\omega,t)] ,$$
$$\widehat{u}_k(\omega,t)=\widehat{\varphi}_k(\omega )e^{
  -\nu t|\omega|^2}+\int_0^t P_k(\omega,\tau)  e^{ -\nu(t-\tau)|\omega|^2 }d\tau .$$  If  
   $ p,p_1,p_2\geq 0,p,p_1+p_2< n,  $ then $$\int_{{\bf{R}}^n}
  \frac{ |\widehat{\varphi}_k(\omega)|}{|\omega|^p} d\omega=\int_{|\omega|\leq 1}
  \frac{ |\widehat{\varphi}_k(\omega)|}{|\omega|^p}d\omega+\int_{|\omega|\geq 1}
  \frac{| \widehat{\varphi}_k(\omega)|}{|\omega|^p}d\omega    \leq \frac{| \varphi|_{0-}S_n  }
  {n-p}+|\varphi|_{0+}.$$ By  H\"{o}lder inequality,we have  \begin{eqnarray*}&\/& \int_{ {{\bf{R}}}^n} 
     \frac{|\widehat{\varphi}_k(\xi)\widehat{\eta}_k(\omega-\xi)|d\xi}{|\xi|^{p_1}
  |\omega-\xi|^{p_2}}\leq
   \left[\int_{ {{\bf{R}}}^n}  \frac{|\widehat{\varphi}_k(\xi)\widehat{\eta}_k(\omega-\xi)|d\xi}{|\xi|^{p_1+p_2} } 
   \right]^{\frac{p_1}{p_1+p_2}} \left[\int_{ {{\bf{R}}}^n}  \frac{|\widehat{\varphi}_k(\xi)
   \widehat{\eta}_k(\omega-\xi)|d\xi}{ |\omega-\xi|^{p_1+p_2}} \right]^{\frac{p_2}{p_1+p_2}}\\  &\leq&
  \left[ \sup\limits_{\xi\in {\bf{R}}^n} |\widehat{\eta}_k(\xi)|\cdot\int_{ {{\bf{R}}}^n}  
  \frac{|\widehat{\varphi}_k(\xi) |d\xi}{|\xi|^{p_1+p_2} } \right]^{\frac{p_1}{p_1+p_2}} \cdot 
  \left[\sup\limits_{\xi\in {\bf{R}}^n} |\widehat{\varphi}_k(\xi)|\cdot\int_{ {{\bf{R}}}^n}  
  \frac{|\widehat{ \eta}_k(\omega-\xi)|d\xi}{ |\omega-\xi|^{p_1+p_2}} \right]^{\frac{p_2}{p_1+p_2}}\\
    &\leq &
  |\eta_k|_{0-} ^{\frac{p_1}{p_1+p_2}}\cdot  | \varphi_k|_{0-}^{\frac{p_2}{p_1+p_2}}\cdot  
  [ \frac{| \varphi_k|_{0-}S_n  }{n-p_1-p_2}+|\varphi_k|_{0+} ]^{\frac{p_1}{p_1+p_2}} \cdot 
  [ \frac{| \eta_k|_{0-}S_n  }{n-p_1-p_2}+|\eta_k|_{0+} ]^{\frac{p_2}{p_1+p_2}}, \\  &\leq &
  ||\eta || _0\cdot  || \varphi|| \cdot (\frac{ S_n  }{n-p_1-p_2}+\lambda ),  \end{eqnarray*}
\begin{eqnarray*}&\/&   (2\pi)^n  \widehat{u_j\frac{\partial \rho_k}
  {\partial  x_j}}(\omega,t)     =
    ( \widehat{u}_j   \ast \widehat{\frac{\partial \rho_k}{\partial  x_j}})(\omega,t) =  \int_{ {{\bf{R}}}^n}
   \widehat{u}_j  (\omega-\xi,t)  \widehat{\frac{\partial \rho_k}{\partial  x_j}}(\xi,t)d\xi
    \\ &=&i
\int_{ {{\bf{R}}}^n}
   \xi_j  \widehat{u}_j  (\omega-\xi,t)   \widehat{\rho}_k  (\xi,t)
     d\xi  \\  &=&i
 \int_{ {{\bf{R}}}^n}   \xi_j       \left[\widehat{\varphi}_j(\omega-\xi)
e^{ -\nu t|\omega-\xi|^2}  +\int_0^t P_j(\omega-\xi,\tau_1)
e^{ -\nu(t-\tau_1)|\omega-\xi|^2 }d\tau_1 \right] \\  &\/&   \left\{\widehat{\eta}_k(\xi )e^{ -\kappa t|
  \xi|^2}+\int_0^t [\widehat{g}_k(\xi,\tau_2)+\widehat{h}_k(\xi,\tau_2)]
   e^{ -\kappa(t-\tau_2)|\xi|^2 }d\tau_2\right\}  d\xi  , \end{eqnarray*} 
     \begin{eqnarray*}&\/&  (2\pi)^n  \widehat{u_j\frac{\partial u_k}
  {\partial  x_j}}(\omega,t) =   \int_{ {{\bf{R}}}^n}
   \widehat{u}_j  (\omega-\xi,t)  \widehat{\frac{\partial u_k}{\partial  x_j}}(\xi,t)d\xi
    \\ &=&i
   \int_{ {{\bf{R}}}^n}
   \xi_j  \widehat{u}_j  (\omega-\xi,t)   \widehat{u}_k  (\xi,t)
     d\xi     \\  &=&i
   \int_{ {{\bf{R}}}^n}   \xi_j        [\widehat{\varphi}_j(\omega-\xi)
e^{ -\nu t|\omega-\xi|^2}   +\int_0^t P_j(\omega-\xi,\tau_1)
e^{ -\nu(t-\tau_1)|\omega-\xi|^2 }d\tau_1   ] \\  &\/&    
[\widehat{\varphi}_k(\xi)e^{ -\nu t|\xi|^2}+\int_0^t P_k(\xi,\tau_2)e^{ -\nu(t-\tau_2)|\xi|^2 }
d\tau _2 ]  d\xi    \end{eqnarray*}  By    $$u_  t + u \cdot
 \nabla u-\nu\triangle u+\nabla p-A
 \rho -f=r+  u \cdot
  \nabla u,$$   $$\rho _t+u\cdot \nabla  \rho -  \kappa  \triangle
   \rho-g=h+  u \cdot  \nabla \rho,$$  As  $\textrm{div}\varphi=0,
   \textrm{div}r(y,t) +
    \textrm{div}f(y,t) +\sum\limits_{i=1}^n\sum\limits_{j=1}^na_{i,j}\frac{\partial  \rho_j(y,t )}{\partial  y_i}=
  \Delta p (y,t)$, then   \begin{eqnarray*}  &\/&\textrm{div}
    u  (x,t)=
  \int_{ {\bf{R}}^n}\frac{e^{{\frac{-|x-y|^2}{4\nu t}}}  \textrm{div}\varphi(y )
   }{(4\pi\nu t)^{ \frac{n}{2}}} dy\\   &\/&+\int_0^t\int_{ {\bf{R}}^n}
    \frac{e^{\frac{-|x-y|^2}{4 \nu(t-\tau)}}[\textrm{div}r(y,\tau)+\textrm{div}f(y,\tau)+\sum\limits_{i=1}^n
    \sum\limits_{j=1}^na_{i,j}\frac{\partial  \rho_j(y,\tau )}{\partial  y_i}-
  \Delta p (y,\tau) ]}{[4\pi\nu (t-\tau)]^{ \frac{n}{2}}}d\tau dy   \\
   &=&0.     \end{eqnarray*}So   $(1)\Longleftrightarrow  r+u\cdot \bigtriangledown  
   u=0,h+u\cdot \bigtriangledown  \rho=0,(2)\Longleftrightarrow  r+u\cdot \bigtriangledown  u=0.$

(I)  (i)  If   $r,h\in   H^n(\lambda,D),$
 \begin{eqnarray*}&\/&   (2\pi)^n   \left|  u_j\frac{\partial \rho_k}
  {\partial  x_j}    \right|_{0-}     \leq
 \int_{ {{\bf{R}}}^n} [| \xi|    \cdot   |\widehat{\eta}_k(\xi )|+  
  \sup\limits_{t\geq 0}\frac{|\widehat{g}_k(\xi,t)|+|\widehat{h}_k(\xi,t)|}{\kappa|\xi|} ] \\&\/&  
   [|\widehat{\varphi}_j(\omega-\xi)|+
 \nu^{-1}(n+1)\sup\limits_{t\geq 0}\frac{| \widehat{r}_j(\omega-\xi,t)|+
 | \widehat{f}_j(\omega-\xi,t)|+nB| \widehat{\eta}_j(\omega-\xi)|}{|\omega-\xi|^2}\\ 
  &\/&+n (n+1) B\nu^{-1}\kappa^{-1}\sup\limits_{t\geq 0,1\leq j\leq n}\frac{|g_j(\omega-\xi,t)|+
  |h_j(\omega-\xi,t)|}{|\omega-\xi|^4}
 ]d\xi \\&\leq&
| \eta|_{1-}[ | \varphi|_{0+} +\nu^{-1}(n+1)(D+||f||+nB||\eta||)(\frac{ S_n  }{ n-2 }+\lambda)
\\&\/&+(n+1)nB\nu^{-1}\kappa^{-1}(D+||g|| )( \frac{S_n  }{n-4} +\lambda)]+  \kappa^{-1} |\varphi|_{0-}(D+||g|| )
(  \frac{S_n  }{n-1} +\lambda) \\  &\/& +  
\nu^{-1}\kappa^{-1} (n+1) (D+||f||  +nB|| \eta||)(D+||g|| )(\frac{S_n  }{n-3} +\lambda)  \\&\/&
 +n (n+1) B\kappa^{-2}\nu^{-1}
   ( ||g|| +D)^2\cdot  (\frac{   S_n}{n-5}+\lambda) , \end{eqnarray*}
\begin{eqnarray*}&\/&   (2\pi)^n \left|   u_j\frac{\partial u_k}
  {\partial  x_j}\right|_{0-}   \leq
 \int_{ {{\bf{R}}}^n}   [|\widehat{\varphi}_k( \xi)|\cdot  |\xi|+
 \nu^{-1}(n+1)\sup\limits_{t\geq 0}\frac{| \widehat{r}_k( \xi,t)|+
 | \widehat{f}_k( \xi,t)|+nB| \widehat{\eta}_k( \xi)|}{| \xi|}   
 \\&\/&+ n(n+1)B\nu^{-1}\kappa^{-1}\sup\limits_{t\geq 0}\frac{|g_k( \xi,t)|+|h_k( \xi,t)|}{| \xi|^3}   ]
    [|\widehat{\varphi}_j(\omega- \xi)|\\&\/&+
 (n+1)\nu^{-1}\sup\limits_{t\geq 0}\frac{| \widehat{r}_j(\omega-\xi,t)|+| \widehat{f}_j(\omega-\xi,t)|+nB| 
 \widehat{\eta}_j(\omega-\xi)|}{|\omega-\xi|^2}\\&\/&  + n(n+1)B\nu^{-1}\kappa^{-1}\sup\limits_{t\geq 0}
 \frac{|g_j(\omega-\xi,t)|+|h_j(\omega-\xi,t)|}{|\omega-\xi|^4}   ]d\xi \\&\leq& | \varphi|_{1-}[|\varphi|_{0+} +(n+1)\nu^{-1}
(D+||f|| +nB|| \eta||)( \frac{S_n  }{ n-2 } +\lambda)\\&\/&+n(n+1)B\nu^{-1}\kappa^{-1}(D+||g|| )
( \frac{S_n  }{n-4} +\lambda)]+(n+1)\nu^{-1}| \varphi|_{0-} (D+||f|| +nB|| \eta||)(\frac{S_n  }{n-1}+\lambda) 
\\&\/&+(n+1)^2\nu^{-2}(||f||   +D+nB| |\eta||_{0})^2 (\frac{  S_n}{n-3}+ \lambda) \\ 
 &\/& + 2 n(n+1)^2B \nu^{-2}\kappa^{-1}  (D+||f||_{0}    +nB ||\eta||_{0})   (D+||g||_{0} )
 (\frac{  S_n }{n-5}+ \lambda)\\&\/&+ n(n+1)|\varphi|_{0-}B\nu^{-1}\kappa^{-1}(D+||g|| )
 ( \frac{S_n  }{n-1} +\lambda) ]\\&\/& + n^2 (n+1)^2B^2\nu^{-2}\kappa^{-2}
  (||g|| +D)^2  (\frac{  S_n}{n-7}+\lambda )  ,\end{eqnarray*}  \begin{eqnarray*}&\/&(2\pi)^n\left|u_j  \frac{\partial
    \rho_k}{\partial  x_j}\right|_{0+} \leq
(    | \varphi _j|_{0+}+ \nu^{-1} |P_j|_{-2+} ) [  |\eta_k|_{1+} + \kappa^{-1}
 (|g_k|_{-1+}+|h_k|_{-1+} )]\\ &\leq& [ |\eta|_{1+} + \kappa^{-1}(||g||  +D)( \frac{ S_n}{ n-1 }+\lambda)]
[ | \varphi|_{0+}  + (n+1)\nu^{-1}(D+ ||f|| + nB||\eta|| ) (  \frac{  S_n}{n-2}+\lambda)\\&\/&+
n(n+1)B\kappa^{-1}\nu^{-1}(||g||+D)(\frac{S_n}{n-4}+\lambda)] ,  \end{eqnarray*}
 \begin{eqnarray*}&\/&(2\pi)^n\left|u_j  \frac{\partial
    u_k}{\partial  x_j}\right|_{0+} \leq
[    | \varphi _j|_{0+} + \nu^{-1} |P_j|_{-2+}] [  |\varphi_k|_{1+} +  \nu^{-1} |P_k|_{-1+}] 
\\ &\leq& [ | \varphi|_{0+}  + (n+1)\nu^{-1}(D+ ||f|| + nB||\eta|| ) (  \frac{  S_n}{n-2}+\lambda)+
n(n+1)B\kappa^{-1}\nu^{-1}(||g||+D)(\frac{S_n}{n-4}+\lambda)]\\&\/&
 [ | \varphi|_{1+}  + (n+1)\nu^{-1}(D+ ||f|| + nB||\eta|| ) (  \frac{  S_n}{n-1}+\lambda) +
n(n+1)B\kappa^{-1}\nu^{-1}(||g||+D)(\frac{S_n}{n-3}+\lambda)]  , 
\end{eqnarray*}$$|(u \cdot \nabla ) u|_{0-}   \leq  (2\pi)^{n} D,|(u 
\cdot \nabla ) \rho|_{0-}   \leq  (2\pi)^{n} D,|(u \cdot \nabla ) u|_{0+}   \leq  (2\pi)^{n} 
\lambda^{-1} D,|(u \cdot \nabla  )\rho|_{0+}   \leq  (2\pi)^{n} \lambda^{-1}   D.$$  If 
  $n\geq 8,1\leq k\leq n,h_k,r_k\in H(\lambda,D),$  then     $ -(u\cdot\nabla ) u,-(u\cdot\nabla )
   \rho  \in  H^{n}(\lambda,D) $.

  (ii)Let   $(r^{[1]},h^{[1]})\in H^{2n}(\lambda,D), (r^{[m+1]},h^{[m+1]})=\{-[u(r^{[m]},h^{[m]})\cdot 
   \bigtriangledown ] u(r^{[m]},h^{[m]}),-[u(r^{[m]},h^{[m]})\cdot  \bigtriangledown ] \rho(r^{[m]},h^{[m]})] \}. $
 Then
  $$ \widehat{\rho}^{[m]}_{k}(\omega,t)-\widehat{\rho}^{[m-1]}_{k}(\omega,t)=
   \int_0^t [\widehat{h}^{[m-1]}_{k}(\omega,\tau)-\widehat{h}^{[m-2]}_{k}(\omega,\tau)]
   e^{ -\kappa(t-\tau)|\omega|^2 }d\tau ,$$   \begin{eqnarray*}&\/& \widehat{u}^{[m]}_{k}(\omega,t)-
    \widehat{u}^{[m-1]}_{k}(\omega,t)=- \sum\limits_{l=1}^n\frac{\omega_l\omega_k}{|\omega|^2} \int_0^t
     \{\widehat{r}^{[m-1]}_{l}(\omega,\tau)-\widehat{r}^{[m-2]}_{l}(\omega,\tau)\\   &\/&  +e_lA
    \int_0^\tau [ \widehat{h}^{[m-1]}_{l}(\omega,\theta)-\widehat{h}^{[m-2]}_{l}(\omega,\theta)]
   e^{ -\kappa(\tau-\theta)|\omega|^2 }d\theta\}
   e^{ -\nu(t-\tau)|\omega|^2 }d\tau  \\   &\/&+ \int_0^t \{\widehat{r}^{[m-1]}_{k}(\omega,\tau)-
   \widehat{r}^{[m-2]}_{k}(\omega,\tau)+e_kA
    \int_0^\tau [ \widehat{h}^{[m-1]}_{k}(\omega,\theta)-\widehat{h}^{[m-2]}_{k}(\omega,\theta)]
   e^{ -\kappa(\tau-\theta)|\omega|^2 }d\theta\}
   e^{ -\nu(t-\tau)|\omega|^2 }d\tau ,  \end{eqnarray*}

 \begin{eqnarray*}&\/&   (2\pi)^n   \left|  ( u_j^{[m]}-u_j^{[m-1]})   \frac{\partial \rho_k^{[m]}}
  {\partial  x_j}       \right|_{0-}   \\ &\leq&   (n+1)\nu^{-1} [
| \eta|_{1-}      (\frac{ S_n  }{ n-2 }+\lambda)+\kappa^{-1}(D+||g||)  
 (\frac{S_n  }{n-3} +\lambda)]||r^{[m-1]}- r^{[m-2]}  || \\&\/&+ (n+1)nB\nu^{-1}\kappa^{-1}[| \eta|_{1-}
 ( \frac{S_n  }{n-4} +\lambda)+  \kappa^{-1}( ||g||+D)(  \frac{S_n  }{n-5} +\lambda) ] ||h^{[m-1]}- h^{[m-2]}
   || , \end{eqnarray*}
    \begin{eqnarray*}&\/&   (2\pi)^n   \left|  u_j^{[m-1]}\frac{\partial (\rho_k^{[m]}-\rho_k^{[m-1]})}
  {\partial  x_j}        \right|_{0-}   \\ &\leq&
  [\kappa^{-1} |\varphi|_{0-}(  \frac{S_n  }{n-1} +\lambda) +  \kappa^{-1}  \nu^{-1}(n+1) (D+||f||  
  +nB|| \eta||)(\frac{S_n  }{n-3} +\lambda)  \\&\/&
 +n (n+1) B\kappa^{-2}\nu^{-1}
   ( ||g|| +D)\cdot  (\frac{   S_n}{n-5}+\lambda)] ||h^{[m-1]}- h^{[m-2]}  || , \end{eqnarray*}
\begin{eqnarray*}&\/&   (2\pi)^n \left| ( u_j^{[m]}- u_j^{[m-1]})\frac{\partial u_k^{[m]}}
  {\partial  x_j}   \right|_{0-}   \\&\leq&  (n+1)[| \varphi|_{1-} \nu^{-1}(\frac{S_n  }{n-2}+\lambda)+
   (n+1)\nu^{-2}\kappa^{-1}(D+||g||)( \frac{S_n  }{n-5} +\lambda) \\&\/&+ (n+1) \nu^{-2} (D+||f||    +nB ||\eta||)
 (\frac{  S_n }{n-3}+ \lambda)]  ||r^{[m-1]}- r^{[m-2]}  ||\\&\/& + (n+1)\kappa  ^{-1}\nu^{-1}[| \varphi|_{1-} 
 (\frac{S_n  }{n-4}+\lambda)+ (n+1)\nu^{-1}\kappa^{-1}(D+||g||)( \frac{S_n  }{n-7} +\lambda) \\&\/&+ (n+1) \nu^{-1} (D+||f||
     +nB ||\eta||)
 (\frac{  S_n }{n-5}+ \lambda)] ||h^{[m-1]}- h^{[m-2]}  ||   ,\end{eqnarray*}
    \begin{eqnarray*}&\/&   (2\pi)^n \left|    u_j^{[m-1]}\frac{\partial ( u_k^{[m]}-u_k^{[m-1]})}
  {\partial  x_j}  \right|_{0-}    \\&\leq& (n+1)\nu^{-1}[| \varphi|_{1-} (\frac{S_n  }{n-2}+\lambda)+ 
  (n+1)\nu^{-1}\kappa^{-1}(D+||g||)( \frac{S_n  }{n-5} +\lambda) \\&\/&+ (n+1) \nu^{-1} (D+||f||    +nB ||\eta||)
 (\frac{  S_n }{n-3}+ \lambda)]  ||r^{[m-1]}- r^{[m-2]}  ||\\&\/& + (n+1)\nu^{-1}\kappa  ^{-1}[
  | \varphi|_{1-} (\frac{S_n  }{n-4}+\lambda)+ (n+1)\nu^{-1}\kappa^{-1}(D+||g||)( \frac{S_n  }{n-7} +\lambda) 
  \\&\/&+ (n+1) \nu^{-1} (D+||f||    +nB ||\eta||)
 (\frac{  S_n }{n-5}+ \lambda)] ||h^{[m-1]}- h^{[m-2]}  ||  ,\end{eqnarray*}
   \begin{eqnarray*}&\/&   (2\pi)^n  | (u_j^{[m]}- u_j^{[m-1]}  )\frac{\partial \rho_k^{[m]}}
  {\partial  x_j} |_{0+} \leq(n+1)\nu^{-2}[ |\eta|_{1+}(\frac{  S_n }{n-2}+ \lambda) + \kappa^{-1}(||g||  +D)
  ( \frac{ S_n}{ n-1 }+\lambda)]\\&\/&[( \frac{ S_n}{ n-2 }+\lambda)||r^{[m-1]}- r^{[m-2]}  ||+\kappa^{-1}(
     \frac{ S_n}{ n-4 }+\lambda)||h^{[m-1]}- h^{[m-2]}  ||]  ,  \end{eqnarray*}
  \begin{eqnarray*}&\/&   (2\pi)^n  |  u_j^{[m-1]} \frac{\partial (\rho_k^{[m]}-\rho_k^{[m-1]})}
  {\partial  x_j} |_{0+} \\  &   \leq&\kappa^{-1}
[ | \varphi|_{0+}  + (n+1)\nu^{-1}(D+ ||f|| + nB||\eta|| ) (  \frac{  S_n}{n-2}+\lambda)+\kappa^{-1}\nu^{-1}
n(n+1)B(||g||+D)(\frac{S_n}{n-4}+\lambda)]\\&\/&(\frac{S_n}{n-1}+\lambda)||h^{[m-1]}- h^{[m-2]}  ||  ,  \end{eqnarray*}
 \begin{eqnarray*}&\/&(2\pi)^n\left|(u_j^{[m]}- u_j^{[m-1]}  )\frac{\partial u_k^{[m]}}
  {\partial  x_j}\right|_{0+} \\ &\leq&  [ | \varphi|_{1+}  + (n+1)\nu^{-1}(D+ ||f|| + nB||\eta|| ) 
  (  \frac{  S_n}{n-1}+\lambda) +\kappa^{-1}\nu^{-1}
n(n+1)B(||g||+D)(\frac{S_n}{n-3}+\lambda)]\\&\/&(n+1)\nu^{-1}
[(\frac{S_n}{n-2}+\lambda)||r^{[m-1]}- r^{[m-2]}  ||+
\kappa^{-1}nB(\frac{S_n}{n-4}+\lambda)||h^{[m-1]}- h^{[m-2]}  || ]  . 
\end{eqnarray*}\begin{eqnarray*}&\/&(2\pi)^n\left| u_j^{[m-1]} \frac{\partial (u_k^{[m]}- u_k^{[m-1]}  )}
  {\partial  x_j}\right|_{0+} \\ &\leq& [ | \varphi|_{0+}  + (n+1)\nu^{-1}(D+ ||f|| + nB||\eta|| ) 
  (  \frac{  S_n}{n-2}+\lambda)+\kappa^{-1}\nu^{-1}
n(n+1)B(||g||+D)(\frac{S_n}{n-4}+\lambda)]\\&\/&(n+1)\nu^{-1}[(\frac{S_n}{n-3}+\lambda)||r^{[m]}- r^{[m-1]}  ||+
\kappa^{-1}nB(\frac{S_n}{n-1}+\lambda)||h^{[m]}- h^{[m-1]}  ||]  . \end{eqnarray*}(iii)\begin{eqnarray*}&\/&
  ||(r^{[m+1]},h^{[m+1]})-(r^{[m]},h^{[m]})||\\&=& ||(-u ^{[m]} \cdot  \bigtriangledown u ^{[m]} ,-u  ^{[m]} 
  \cdot  \bigtriangledown  \rho ^{[m]} ) +(u ^{[m-1]} \cdot  \bigtriangledown u^{[m-1]},u  ^{[m-1]} \cdot  
  \bigtriangledown  \rho ^{[m-1]} )|| \\  &\leq &||[ (u ^{[m]} -u ^{[m-1]})\cdot   \bigtriangledown u ^{[m]}
     ,(u ^{[m]} -u ^{[m-1]})\cdot   \bigtriangledown \rho ^{[m]}]||\\&\/&+ ||[  u ^{[m-1]}  \cdot 
       \bigtriangledown (u ^{[m]} - u ^{[m-1]}) ,u ^{[m-1]}  \cdot   \bigtriangledown (\rho ^{[m]} - 
       \rho ^{[m-1]})]|| \\&\leq &   \varepsilon   ||(r^{[m]},h^{[m]})-(r^{[m-1]},h^{[m-1]})||, 
      \end{eqnarray*}By compression  fixed point principle
   there  exist  solutions   $p\in  C^\infty({\bf{R}}^n), u_k,\rho_k\in H(1\leq k\leq n) $   for system (1)   as $n\geq 8 $.

     (II)(i)If      $\rho(x,t)=0,g(x)=\eta(x)=0,   $  then
 \begin{eqnarray*}&\/&\widehat{u}_k(\omega,t)=\widehat{\varphi}_k(\omega )e^{
  -\nu t|\omega|^2}+\int_0^t [\widehat{r}_k(\omega,\tau)+\widehat{f}_k(\omega,\tau)
-\sum\limits_{m=1}^n\frac{\omega_m\omega_k}{|\omega|^2}
        [  \widehat{r_m}(\omega,\tau)+ \widehat{f_m}(\omega,\tau)]]
e^{ -\nu(t-\tau)|\omega|^2 }d\tau  ,\end{eqnarray*}
  \begin{eqnarray*}&\/& (2\pi)^n  \widehat{u_j\frac{\partial u_k}
  {\partial  x_j}} (\omega,t) =   \int_{ {{\bf{R}}}^n}
   \widehat{u}_j  (\omega-\xi,t)  \widehat{\frac{\partial u_k}{\partial  x_j}}(\xi,t)d\xi
    \\ &=&
   \int_{ {{\bf{R}}}^n}
   \xi_j  \widehat{u}_j  (\omega-\xi,t)   \widehat{u}_k  (\xi,t)
     d\xi     \\  &=&
   \int_{ {{\bf{R}}}^n}   \xi_j        \{\widehat{\varphi}_j(\omega-\xi)
e^{ -\nu t|\omega-\xi|^2}   +\int_0^t [\widehat{r}_j(\omega-\xi,\tau_1)+\widehat{f}_j(\omega-\xi,\tau_1)
]
e^{ -\nu(t-\tau_1)|\omega-\xi|^2 }d\tau_1  \} \\  &\/&  
  \{\widehat{\varphi}_k(\xi)e^{ -\nu t|\xi|^2}+\int_0^t [\widehat{r}_k( \xi,\tau_2)+\widehat{f}_k( \xi,\tau_2)
]  e^{ -\nu(t-\tau_2)|\xi|^2 }
d\tau _2 \}  d\xi    \end{eqnarray*}
\begin{eqnarray*}&\/&  (2\pi)^n  \left|  u_j\frac{\partial u_k}
  {\partial  x_j}  \right|_{0-}  \\&\leq&
 \int_{ {{\bf{R}}}^n}   [|\widehat{\varphi}_k( \xi)|\cdot  |\xi|+
 \nu^{-1}\sup\limits_{t\in  {\mathbf{R}}_+}\frac{|\widehat{r}_k( \xi,t)|+|\widehat{f}_k( \xi,t)|}{| \xi|}    ] 
 \\  &\/&    [|\widehat{\varphi}_j(\omega-\xi)|+
 \nu^{-1}\sup\limits_{}\frac{|\widehat{r}_j(\omega- \xi,t)|+|\widehat{f}_j( \omega- \xi,t)| }{|\omega-\xi|^2}  ]d\xi \\&\leq&
| \varphi|_{1-}|\varphi|_{0+}+\nu^{-1}| \varphi|_{1-}(D+||f||  )(  \frac{  S_n}{ n-2 }+  \lambda) 
\\ &\/& +\nu^{-1}| \varphi|_{0-}(D+||f||  )(  \frac{  S_n}{ n-1 }+  \lambda) +\nu^{-2}(  \frac{  S_n}{ n-3 }+
  \lambda)(||f||   +D)^2 \leq \frac{(2\pi)^n D}{n},\end{eqnarray*}
 \begin{eqnarray*}&\/&(2\pi)^n\left|u_j  \frac{\partial
    u_k}{\partial  x_j}\right|_{0+} \\ &\leq&
[ | \varphi|_{0+}  + \nu^{-1}(D+||f||  )(  \frac{  S_n}{ n-2 }+  \lambda) ][ |\varphi|_{1+}  + \nu^{-1}(D+||f||  )(
    \frac{  S_n}{ n-1 }+  \lambda)]\leq \frac{\lambda^{-1}(2\pi)^n   D}{n}, \end{eqnarray*}   
    So  if  $r\in  H^n(\lambda,D),n\geq  4,$   $-(u\cdot\nabla)  u\in H^n(\lambda,D).$

(ii)If  $r^{[1]}\in  H^n(\lambda,D),n\geq  4,$   then
\begin{eqnarray*}&\/&  (2\pi)^n  \left| (  u_j^{[m]}- u_j^{[m-1]}
 )\frac{\partial u_k^{[m]}}
  {\partial  x_j}  \right|_{0-}  \\&\leq&
 \int_{ {{\bf{R}}}^n}
  \sup\limits_{t\geq 0}\frac{|\widehat{r}_k^{[m-1]}( \xi,t)- \widehat{r}_k^{[m-2]}( \xi,t)|}{\nu| \xi|^2}    
  [|\widehat{\varphi}_j(\omega-\xi)||\omega-\xi|+
 \nu^{-1}\sup\limits_{t\geq 0}\frac{|\widehat{r}_j(\omega- \xi,t)|+|\widehat{f}_j( \omega- \xi,t)| }{|\omega-\xi|}  
 ]d\xi \\&\leq& [\nu^{-1}| \varphi|_{1-}(  \frac{  S_n}{ n-2 }+  \lambda) +\nu^{-2}(D+||f||  )(  \frac{  S_n}{ n-3 }+ 
  \lambda)]||r^{[m-1]}-r^{[m-2]}||,\end{eqnarray*}

  \begin{eqnarray*}&\/& (2\pi)^n \left|  u_j(\frac{\partial u_k^{[m]} }
  {\partial  x_j}- \frac{\partial u_k^{[m-1]}}
  {\partial  x_j})   \right|_{0-} \\&\leq&
 \int_{ {{\bf{R}}}^n}
  \sup\limits_{t\geq 0}\frac{|\widehat{r}_k^{[m-1]}( \xi,t)- \widehat{r}_k^{[m-2]}( \xi,t)|}{\nu| \xi|}   
   [|\widehat{\varphi}_j(\omega-\xi)|+
 \nu^{-1}\sup\limits_{t\geq 0}\frac{|\widehat{r}_j(\omega- \xi,t)|+|\widehat{f}_j( \omega- \xi,t)| }{|\omega-\xi|^2} 
  ]d\xi  \\ &\leq&
 \nu^{-1}[(  \frac{  S_n}{ n-1 }+  \lambda)| \varphi|_{0-} +\nu^{-1}(D+||f||  )(  \frac{  S_n}{ n-3 }+  
 \lambda)]||r^{[m-1]}-r^{[m-2]}||\end{eqnarray*}

 \begin{eqnarray*}&\/&  (2\pi)^n  \left| (  u_j^{[m]}- u_j^{[m-1]}
 )\frac{\partial u_k^{[m]}}
  {\partial  x_j}  \right|_{0+}  \\&\leq& \nu^{-1}(  \frac{  S_n}{ n-2 }+  \lambda)[ | \varphi|_{1+}+\nu^{-1}(
      \frac{  S_n}{ n-1 }+  \lambda) (D+||f||) ]||r^{[m-1]}-r^{[m-2]}||,\end{eqnarray*}

 \begin{eqnarray*}&\/& (2\pi)^n \left|  u_j(\frac{\partial u_k^{[m]} }
  {\partial  x_j}- \frac{\partial u_k^{[m-1]}}
  {\partial  x_j})   \right|_{0+}  \\ &\leq&
 \nu^{-1}[(  \frac{  S_n}{ n-1 }+  \lambda)| \varphi|_{0+} +\nu^{-1}(D+||f||  )(  \frac{  S_n}{ n-3 }+  
 \lambda)]||r^{[m-1]}-r^{[m-2]}||\end{eqnarray*}

(iii)\begin{eqnarray*}&\/&|| r^{[m+1]} - r^{[m]} ||= || -(u ^{[m]} \cdot  \bigtriangledown) u ^{[m]}  
   +(u ^{[m-1]} \cdot  \bigtriangledown )u^{[m-1]} || \\  &\leq &||[ (u ^{[m]} -u ^{[m-1]})\cdot 
     \bigtriangledown ]u ^{[m]}  || + ||[  u ^{[m-1]}  \cdot   \bigtriangledown ](u ^{[m]} - u ^{[m-1]})  || 
      \\&\leq &   \delta   || r^{[m]} - r^{[m-1]} || , \end{eqnarray*}

    By compression  fixed point principle
   there  exist  solutions   $p\in  C^\infty({\bf{R}}^n), u_k,\rho_k\in H(1\leq k\leq n) $   for system (2)   as $n\geq 4 $.
\\

\begin{theorem}
 Let   $   1\leq k\leq n,
 \nu,C>0,f_k, g_k\in K,\varphi_k, \eta_k\in J,
    \varphi=(\varphi_1,\cdots, \varphi_n),\textrm{div} \varphi=0, \overline{I}(\varphi)= 
    \sum\limits_{ m\in{\bf{Z}}^n } |\widehat{\varphi}(m)| , I(\psi)= \sum\limits_{ m\in{\bf{Z}}^n } 
    \int_0^\infty |\widehat{\psi}(m,t)|dt,K(C)= K\cap  \{\alpha| I(\alpha)\leq C\},J(C)= J\cap  
    \{\alpha| \overline{I}(\alpha)\leq C\},$
$$K=\{\alpha\in  C^{\infty}({\mathbf{R}^n\times{\mathbf{R}^+}})|\alpha (x,t)=\sum\limits_{\beta\in 
 {\mathbf{Z}^n}} \widehat{\alpha}(\beta,t)\exp(\sum\limits_{j=1}^n\frac{2\beta_j\pi  x_ji}{l_j}),I(\alpha)<+\infty\},$$
 $$J=\{\alpha\in  C^{\infty}({\mathbf{R}^n})|\alpha (x)=\sum\limits_{\beta\in  {\mathbf{Z}^n}} 
 \widehat{\alpha}(\beta)\exp(\sum\limits_{j=1}^n\frac{2\beta_j\pi  x_ji}{l_j}),\overline{I}(\alpha)<+\infty\},$$
$$
P=\{\alpha|\alpha(x,t)= \sum\limits_{\beta\in{\bf{Z}}^n} \widehat{\alpha}(\beta,t)\exp(
 \sum\limits_{j=1}^n  \frac{2i\beta_j\pi x_j}{l_j} ),
 \alpha\in  C^{\infty}({\mathbf{R}^n\times{\mathbf{R}^+}})\},$$
$$ \varepsilon:= \frac{\sum\limits_{k=1}^nl_k}{\pi \nu} [ \overline{I}(\varphi)+(n+1)
I(f)
  +(n+1)
C]<1,\hspace{0.3cm}   \frac{\sum\limits_{k=1}^nl_k}{2\pi \nu} [ \overline{I}(\varphi)+(n+1)
I(f)
  +(n+1)
C]^
2\leq  C.$$
    Then there  exist  solutions   $ u_k\in  K(1\leq k\leq n),p\in P $   for system (2)   as $n\geq 3  $ .\end{theorem}

 Proof:\hspace{0.3cm}  (i)Let   $ 1\leq k\leq n,r=( r_1,\cdots,
               r_n)^\top\in K^n(C), ||r||=  \max\limits_{ 1\leq k\leq n}  I(r_k),$   
                $$r_k(x,t) = \sum\limits_{\beta\in{\bf{Z}}^n}   \widehat{r} _k (\beta,t
) \exp(
 \sum\limits_{j=1}^n   \frac{2i\beta_j\pi x_j}{l_j} ),  \varphi_k(x) = \sum\limits_{\beta\in{\bf{Z}}^n}  
  \widehat{\varphi} _k (\beta
) \exp(
 \sum\limits_{j=1}^n   \frac{2i\beta_j\pi x_j}{l_j} ),$$
\begin{eqnarray*} &\/&u_{k}(x,t)=   \sum\limits_{m\in{\bf{Z}}^n}
  \exp [\sum\limits_{j=1}^n (  \frac{2im_j\pi x_j}{l_j} -
 (\frac{2m_j\pi}{l_j})^2\nu t)]   \{ \widehat{ \varphi}_{k}({\mathbf{m}})
  + \int_0^t[\widehat{r}_{k}({\mathbf{m}},\tau
)+\widehat{f}_{k}({\mathbf{m}},\tau
)\\  &\/&-\widehat{\frac{ \partial p}{\partial x_k}}({\mathbf{m}},t
)]  \exp [\sum\limits_
 {j=1}^n
 (\frac{2m_j\pi}{l_j})^2\nu \tau] d\tau   \}\\&=&   \sum\limits_{m\in{\bf{Z}}^n}
  \exp [\sum\limits_{j=1}^n (  \frac{2im_j\pi x_j}{l_j} -
 (\frac{2m_j\pi}{l_j})^2\nu t)]   \{ \widehat{ \varphi}_{k}({\mathbf{m}})
  + \int_0^t[\widehat{r}_{k}({\mathbf{m}},\tau
)+\widehat{f}_{k}({\mathbf{m}},\tau
)\\  &\/&+
\frac{ im_k  [  \widehat{\textrm{div}r}({\mathbf{m}},t
) + \widehat{\textrm{div}f}({\mathbf{m}},t
)]      }{2\pi l_k\sum\limits_{j=1}^n
   \frac{m_j^2}{l_j^2} }]  \exp [\sum\limits_
 {j=1}^n
 (\frac{2m_j\pi}{l_j})^2\nu \tau] d\tau   \}.
  \end{eqnarray*}
\begin{eqnarray*} &\/&p(x,t)= p_{0}(t) \\   &\/&-\sum\limits_{m\in{\bf{Z}}^n/\{\mathbf{0}\}   }
\frac{ [  \widehat{\textrm{div}r}({\mathbf{m}},t
) + \widehat{\textrm{div}f}({\mathbf{m}},t
) ][\exp({\sum\limits_{j=1}^n  \frac{2im_j\pi x_j }{l_j}})-
\exp({\sum\limits_{j=1}^n  \frac{2im_j\pi x_{0,j} }{l_j}})]  }{4\pi^2\sum\limits_{j=1}^n
   \frac{m_j^2}{l_j^2} }    , \end{eqnarray*}  then  \begin{eqnarray*}\widehat{\frac{ 
    \partial p}{\partial x_k}}({\mathbf{m}},t
)=\left\{ \begin{array}{ccccc}0&,      &
    m=\mathbf{0},\\  -
\frac{ im_k  [  \widehat{\textrm{div}}r({\mathbf{m}},t) + \widehat{\textrm{div} }f({\mathbf{m}},t) 
  }{2\pi l_k\sum\limits_{j=1}^n
   \frac{m_j^2}{l_j^2} } &,
        &  m\neq   \mathbf{0}. \end{array} \right.\end{eqnarray*}

(ii) As  $\textrm{div}\varphi=0,
   \textrm{div}r(y,t) +
    \textrm{div}f(y,t)  =
  \Delta p (y,t)$, so $$  \textrm{ div} \varphi(x) = \sum\limits_{\beta\in{\bf{Z}}^n}  
   \widehat{ \textrm{ div} \varphi(\beta)} \exp(
 \sum\limits_{j=1}^n   \frac{2i\beta_j\pi x_j}{l_j} )=\sum\limits_{\beta\in{\bf{Z}}^n}  
  \sum\limits_{k=1}^n   \frac{2i\beta_k\pi }{l_k}  \widehat{\varphi} _k (\beta
) \exp(
 \sum\limits_{j=1}^n   \frac{2i\beta_j\pi x_j}{l_j} ) , $$$$  \widehat{ \textrm{ div} 
 \varphi(\beta)} =\sum\limits_{k=1}^n   \frac{2i\beta_k\pi }{l_k}  \widehat{\varphi} _k (\beta
),$$ \begin{eqnarray*} &\/&\textrm{div}
    u  (x,t) =   \sum\limits_{m\in{\bf{Z}}^n}
  \exp [\sum\limits_{j=1}^n (  \frac{2im_j\pi x_j}{l_j} -
 (\frac{2m_j\pi}{l_j})^2\nu t)]   \{\sum\limits_{k=1}^n   \frac{2im_k\pi }{l_k} \widehat{ \varphi}_{k}({\mathbf{m}})\\&\/&
  + \int_0^t\sum\limits_{k=1}^n   \frac{2im_k\pi }{l_k}[\widehat{r}_{k}({\mathbf{m}},\tau
)+\widehat{f}_{k}({\mathbf{m}},\tau
)-\widehat{\frac{ \partial p}{\partial x_k}}({\mathbf{m}},t
) \exp [\sum\limits_
 {j=1}^n
 (\frac{2m_j\pi}{l_j})^2\nu \tau] d\tau   \}\\  &=&   \sum\limits_{m\in{\bf{Z}}^n}
  \exp [\sum\limits_{j=1}^n (  \frac{2im_j\pi x_j}{l_j} -
 (\frac{2m_j\pi}{l_j})^2\nu t)]   \{ \widehat{{\textrm{div}}   \varphi } ({\mathbf{m}})\\&\/&
  + \int_0^t[\widehat{{\textrm{div}}   r} ({\mathbf{m}},\tau
)+\widehat{{\textrm{div}}   f }  ({\mathbf{m}},\tau
)-\widehat{\Delta p} ({\mathbf{m}},t
) ]\exp [\sum\limits_
 {j=1}^n
 (\frac{2m_j\pi}{l_j})^2\nu \tau] d\tau   \}.
  \end{eqnarray*} By    $$u_  t + u \cdot
 \nabla u-v\triangle u+\nabla p -f=r+  u \cdot
  \nabla u,$$  $(1) \Longleftrightarrow   r+  u \cdot
  \nabla u=0.$

  (iii)By Lemma 2.1,we have
 \begin{eqnarray*}&\/& ||-u\cdot\nabla  u||  \leq   \frac{\sum\limits_{k=1}^nl_k}{2\pi \nu} [ \overline{I}(\varphi)+ (n+1)
I(f)
  +(n+1)
C]^
2\leq  C.\end{eqnarray*}  So  if  $r\in  K(C),$   then  $-(u\cdot\nabla)  u\in  K(C).$

  (iv)         Let   $r^{[l]}=-    ( u ^{[l-1]}\nabla ) u ^{[l-1]}
,$  \begin{eqnarray*} &\/&u_{k}^{[l]}(x,t)-u_{k}^{[l-1]}(x,t)  =   \sum\limits_{m\in{\bf{Z}}^n}
  \exp [\sum\limits_{j=1}^n (  \frac{2im_j\pi x_j}{l_j} -
 (\frac{2m_j\pi}{l_j})^2\nu t)]   \{ \int_0^t[\widehat{r}_{k}^{ [l-1]}({\mathbf{m}},\tau
)-\widehat{r}_{k}^{ [l-2]}({\mathbf{m}},\tau
)\\  &\/&+
\frac{ im_k  [  \widehat{\textrm{div}r^{ [l-1]}}({\mathbf{m}},t
)-\widehat{\textrm{div}r^{ [l-2]}}({\mathbf{m}},t
) ]      }{2\pi l_k\sum\limits_{j=1}^n
   \frac{m_j^2}{l_j^2} }]  \exp [\sum\limits_
 {j=1}^n
 (\frac{2m_j\pi}{l_j})^2\nu \tau] d\tau   \}.
  \end{eqnarray*} \begin{eqnarray*}&\/&||r^{[l+1]}-r^{[l]}||\\  &\leq &||   
    ( u ^{[l]}\cdot\nabla  )u ^{[l]} -  (u ^{[l]}\cdot\nabla ) u ^{[l-1]} ||+ ||  
     ( u ^{[l]}\cdot\nabla  )u ^{[l-1]} -   ( u ^{[l-1]}\cdot\nabla  )u ^{[l-1]} ||\\  &\leq  &    
     \frac{(n+1)\sum\limits_{k=1}^nl_k}{ \pi \nu} [ \overline{I}(\varphi)+(n+1)
I(f)
  +(n+1)
C]||r^{[l]}-r^{[l-1]}||=\varepsilon||r^{[l]}-r^{[l-1]}||
\end{eqnarray*}  By compression  fixed point principle
   there  exist  solutions  for   $r+  u \cdot
  \nabla u=0$   for   $n\geq 3$.    So  there  exist  $p\in  P, u_k \in K(1\leq k\leq n) $  
   for system (2)   as $n\geq 3 $.\\

    {\bf   Remark  }\hspace{0.3cm}
      If  $f(x,t) \equiv  O,\varphi(x) \equiv  O,$  then   there   exist  $\lambda,D>0$   
      in  Theorem 3.1.  and   $C>0$  in  Theorem 3.2.

\end{document}